\documentstyle{amsppt}
\magnification1200
\NoBlackBoxes
\TagsOnRight
\topmatter
\title{On Non-intersecting Arithmetic Progressions}\endtitle
\author Ernest S. Croot III \endauthor
\abstract
Let $L(c,x) = e^{c\sqrt{\log{x}\log\log{x}}}$.  
We prove that if $a_1 \pmod{q_1},..., a_k \pmod{q_k}$ are a maximal collection 
of non-intersecting arithmetic progressions, with $2 \leq q_1 < q_2 < \cdots
< q_k \leq x$, then 
$$
{x \over L(\sqrt{2}+o(1),x)}\ <\ k\ <\ {x \over L(1/6-o(1),x)}.
$$
In the case for when the $q_i$'s are square-free, we obtain the improved
upper bound
$$
k\ <\ {x \over L(1/2-o(1),x)}.
$$ 
\endabstract
\endtopmatter

\document

\head I.  Introduction \endhead

Suppose that $a_1 \pmod{q_1}, a_2 \pmod{q_2}, ..., a_k \pmod{q_k}$ is a
collection of arithmetic progressions, where $2 \leq q_1 < \cdots < q_k\leq x$, 
with the property that
$$
\{a_i \pmod{q_i}\}\ \cap\ \{a_j \pmod{q_j}\}\  =\ 
\emptyset,\ \text{if $i \neq j$.}
$$
We say that such a collection of arithmetic progressions is disjoint or
non-intersecting.
Let $f(x)$ be the maximum value for $k$, maximized over all choices of 
progressions $a_i \pmod{q_i}$.  Define
$$
L(c,x)\ :=\ \exp(c \sqrt{\log{x}\log\log{x}}),
$$ 
and define 
$$
\align
&\psi(x,y)\ :=\ \#\{ n \leq y\ :\ p \text{ prime},\  p|n\Longrightarrow 
p \leq y\}, \text{ and}\\
&\psi^*(x,y)\ :=\ \#\{n \leq y\ :\ p \text{ prime},\ 
p^a|n\Longrightarrow p^a\leq y\}.
\endalign
$$ 
In [3], Erd\H{o}s and Szemer\'{e}di prove that 
$$
{x \over \exp\left ((\log{x})^{1/2+\epsilon}\right )}\ <\ f(x)\ <\ 
{x \over (\log{x})^c},
$$
for some constant $c > 0$.  (This result is also mentioned in [2]. ) 
Their lower bound can be refined by using more
exact estimates for $\psi(x,L(c,x))$ than was used in their paper.  
Specifically, as direct consequence of [Lemma 3.1, 1], we have the following 
estimate
\proclaim{Lemma 1} For any constant $c > 0$,
$$
\psi(x,L(c,x))\ =\ {x \over L\left (1/(2c)+o(1),x\right)}.\eqno{(1)}
$$ 
\endproclaim
We also have the same estimate for $\psi^*(x,L(c,x))$, since
$$
\eqalign{
\psi(x,L(c,x)) &> \psi^*(x,L(c,x)) > \psi(x,L(c,x))\ -\ \sum_{n^2 > L(c,x)}
\psi(x/n^2,L(c,x)) \cr
&= \psi(x,L(c,x)) - O\left ( {x \over L\left ( 
c/2 + 1/(2c) + o(1), x\right )} \right ),
}\eqno{(2)}
$$ 
Now, let $p$ be the largest prime number less than or equal $L(1/\sqrt{2},x)$.
Let $q_1,q_2,...,q_t$ be the collection of all integers $\leq x$ which are
divisible by $p$, and whose prime power factors are all $< p$.  
From (1) and (2), we have that
$t = x / \left (p L(1/\sqrt{2}+o(1),x) \right ) = x/L(\sqrt{2}+o(1),x)$.
For each $q_i = p \ell_r^{h_r} \ell_{r-1}^{h_{r-1}} \cdots \ell_1^{h_1}$, where 
$p > \ell_r^{h_r} > \ell_{r-1}^{h_{r-1}} > \cdots > \ell_1^{h_1}$ are the
powers of the disctint primes dividing $q_i$,
we choose the residue class $a_i \pmod{q_i}$ using the Chinese Remainder 
Theorem as follows:
$$
\eqalign{
a_i \equiv \ell_r^{h_r} \pmod{p};&\ a_i \equiv \ell_{j-1}^{h_{j-1}} 
\pmod{\ell_j^{h_j}},\ \text{for $2 \leq j \leq r$};\cr
&\text{and finally, }a_i \equiv 0 \pmod{\ell_1^{h_1}}.
}
$$ 
This is exactly the construction which appears in [3] (except that their 
progressions were all square-free), and it is easy to see
that our choice of progressions $a_i \pmod{q_i}$ are disjoint.  Thus, we have
that
$$
f(x) > {x \over L(\sqrt{2}+o(1),x)}.
$$
In this paper we will prove the following results:
\proclaim{Theorem 1}  If $a_1 \pmod{q_1},...,a_k \pmod{q_k}$ are a collection
of disjoint arithmetic progressions, where the $q_i$'s are square-free and
$2 \leq q_1 < \cdots < q_k \leq x$, then
$$
k < {x \over L(1/2 - o(1), x)}.
$$
\endproclaim
\proclaim{Corollary to Theorem 1} 
$$
f(x) < {x \over L(1/6 - o(1),x)}.
$$
\endproclaim
Thus, we will have shown that
$$
{x \over L(\sqrt{2}+o(1),x)}\ <\ f(x)\ <\ {x \over L(1/6-o(1),x)}.
$$
To see how the Corollary follows from Theorem 1, let 
$b_1 \pmod{r_1},...,b_{f(x)} \pmod{r_{f(x)}}$ be a maximal collection of
disjoint arithmetic progressions with $2 \leq r_1 < \cdots < r_{f(x)}
\leq x$.  Suppose, for proof by contradicition, that for some 
$\epsilon < 1/6$
$$
f(x)\ >\ {x \over L\left (1/6 - \epsilon, x \right )}.\eqno{(3)} 
$$
Write each $r_i = \alpha_i\beta_i$, where $\beta_i$ is square-free, 
gcd$(\alpha_i,\beta_i)=1$, and every prime dividing $\alpha_i$ divides to
a power $\geq 2$.  (Note:  we may have $\alpha_i$ or $\beta_i = 1$.)
Now, at least half of $\alpha_i$'s must be $\leq L(1/3,x)$, for
if not we would have from our assumption (3) that
$$
\eqalign{
{x \over 2L(1/6 -\epsilon,x )} < f(x)/2 &<
\#\{r_i\ :\ \alpha_i > L(1/3,x)\}\cr
&<\ x\sum_{n^2 > L(1/3,x)} {1 \over n^2} 
\prod_{\text{$p$ prime}} \left ( 1 + {1 \over p^2} + {1\over p^3} + 
\cdots \right )\cr
&\ll {x \over L(1/6,x)},
} 
$$ 
which is impossible for $x$ large enough in terms of $\epsilon$.  Thus, we
must have that there exists an $\alpha < L(1/3,x)$ for which at least
$f(x)/\left (2L(1/3,x)\right )$ of the $r_i$'s have 
$\alpha_i = \alpha$.
Let $R(\alpha) \subseteq \{r_1,...,r_{f(x)}\}$ be such a collection of 
$r_i$'s, where 
$$
|R(\alpha)| > {f(x)\over 2L(1/3,x)} > 
{x \over 2L(1/2 - \epsilon,x)},
$$
where this last inequality follows from our assumption (3).
Now there must exist a residue class $b \pmod{\alpha}$ for which at least 
$|R(\alpha)|/\alpha$ of the progressions $b_i \pmod{r_i}$ satisfy 
$$
r_i \in R(\alpha),\text{ and } b_i \equiv b \pmod{\alpha}.\eqno{(4)}
$$

Thus, the arithmetic progressions $b_i \pmod{r_i/\alpha}$,
where $r_i$ satisfies (4), is a collection of 
$\geq |R(\alpha)|/\alpha \gg x / (\alpha L(1/2 - \epsilon,x))$ disjoint  
progressions, with distinct square-free moduli $\leq x/\alpha$.  
This contradicts 
Theorem 1 for $x$ sufficiently large in terms of $\epsilon$.  We must
conclude, therefore, that the bound in (3) is false for all 
$\epsilon < 1/6$ and $x > x_0(\epsilon)$, and so the Corollary to 
Theorem 1 follows.    

\head II.  Proof of Theorem 1 \endhead

Before we prove Theorem 1, we will need the following lemma:
\proclaim{Lemma 2}  There are at most $x/L(c/2+o(1),x)$ positive
integers $n \leq x$ such that $\omega(n) > c\sqrt{\log{x}/\log\log{x}}$.
(Recall:  $\omega(n) = \sum_{p | n,\ p\ \text{prime}} 1$.), where $c$ is
some positive constant. 
\endproclaim
\demo{Proof of Lemma 2}  We observe that
$$
\eqalign{
\#\{n \leq x\ :\ \omega(n) > c\sqrt{\log{x}/\log\log{x}}\}\ &<\ 
x\sum_{j > c\sqrt{{\log{x} \over \log\log{x}}}} {\left ( \sum_{p^a \leq x\atop p
\text{ prime}} {1 \over p^a} \right )^j \over j!}\cr
&= {x \over (c \sqrt{\log{x}/\log\log{x}})^
{\{ c+o(1)\}\sqrt{\log{x}/\log\log{x}}}}\cr
&= {x \over L(c/2+o(1),x)}.
}
$$ 

\enddemo

We now resume the proof of Theorem 1.  Consider the
collection of all the $q_i$'s with the properties\medskip

A.  $\omega(q_i) < \sqrt{\log{x} \over \log\log{x}}$, and \smallskip

B.  There exists a prime $p > L(1,x)$,
such that $p|q_i$,\medskip\noindent
Let $\{r_1,...,r_{k'}\}$ be the collection of all such
$q_i$'s satisfying A and B, and where $\{b(r_1),...,b(r_{k'})\}$ are 
their corresponding residue classes.\smallskip

To prove our theorem, we start with the set $S_0 = \{r_1,...,r_{k'}\}$, and
construct a sequence of subsets $S_0 \supseteq S_1 \supseteq S_2 \supseteq 
\cdots$, and a sequence of primes  $p_1,p_2,...$ 
(and let $p_0 = 1$), such that the for each $i \geq 1$, the following
three properties hold\medskip

1.  Each member of $S_i$ is divisible by the primes $p_1,...,p_i$,\smallskip

2.  There exists an integer $A_i$, such that for each $r_j \in S_i$, we have
that $b(r_j) \equiv A_i \pmod{p_1p_2\cdots p_i}$.\smallskip

3.  $|S_i| > |S_{i-1}| / (p_i \sqrt{ \log{x} / \log\log{x}}).$
\medskip\noindent
We continue constructing these subsets until we reach a subset $S_t$ which
has the additional property: \medskip

4.  There exists a prime $p \neq p_1,...,p_t$, $p \geq L(1,x)$
such that at least $|S_t| / \sqrt{ \log{x}/\log\log{x}}$ of the elements
of $S_t$ are divisible by $p$.\medskip\noindent
Let us suppose for the time being that we can construct these sets $S_1,..., 
S_t$.  Applying Property 3 iteratively, together with Property 4, we have
that the number of elements of $S_t$ which are divisible by $p$ (which are
already divisible by $p_1p_2\cdots p_t$ by Property 1) is at least 
$$
{|S_0| \over p_1p_2\cdots p_t (\sqrt{\log{x}/\log\log{x}})^{t+1}} \geq
{|S_0| \over p_1p_2\cdots p_t L(1/2+o(1),x)}, 
$$
(Note: By Property A above we have that $t \leq \sqrt{\log{x}/\log\log{x}}$ 
since every element of $S_0$ has at most $\sqrt{\log{x}/\log\log{x}}$ prime
factors.)
From this, together with the fact that $p > L(1,x)$, we have 
$$
\eqalign{
{x \over p_1\cdots p_k L(1,x)} &\geq 
\#\{n \leq x\ :\ p p_1p_2 \cdots p_t | n\} > \#\{q \in S_t\ :\ p | q\} \cr
&\geq {|S_0| \over p_1p_2\cdots p_t L(1/2+o(1),x)}.
}
$$
It follows that
$$
|S_0| < {x \over L(1/2-o(1),x)},
$$
From this, together with Lemmas 1 and 2 and the fact that the elements of 
$S_0$ satisfy A and B above, we have that
$$
\eqalign{
{x \over L(1/2-o(1),x)} > |S_0| &> k - \#\{n \leq x: 
\omega(n) \geq \sqrt{\log{x}/\log\log{x}}\}\cr
&\ \ \ \ \ \ \ \ \ \ \  - \psi(x,L(1,x))\cr
&> k - {x \over L(1/2-o(1),x)},
}
$$
and so
$$
k < {x \over L(1/2-o(1),x)},
$$
which proves our theorem.\smallskip

To construct our sets $S_i$, we apply the following iterative procedure:  
suppose we have constructed the sets $S_1,...,S_i$, which satisfy 
1 through 3 as above.
To construct $S_{i+1}$, first pick any element $r \in S_i$. 
Now let $e_1,...,e_j$ be all those primes dividing $r/(p_1\cdots p_i)$
(note: $j < \sqrt{\log{x}/\log\log{x}}$).  
Each element $s \in S_i$, $s \neq r$, is divisible by at least one of these
primes, since otherwise 
$\text{gcd}(r,s)=p_1\cdots p_i$ and so we would have 
$b(r) \equiv A_i \equiv b(s) \pmod{\text{gcd}(r,s)}$, which would mean that
$\{ b(r) \pmod{r} \} \cap \{b(s) \pmod{s}\} \neq \emptyset$.\smallskip

Now, there must be at least $|S_i|/j > |S_i|/\sqrt{\log{x}/\log\log{x}}$ 
of the elements of $S_i$ which are divisible by one of these primes $e_h$.  
Let $C_i \subseteq S_i$ be the collection of
all elements $S_i$ divisible by this prime $e_h$.
There exists at least one residue class $B \pmod{e_h}$ for which more than
$|C_i|/e_h > |S_i|/(e_h \sqrt{\log{x}/\log\log{x}})$ of the elements 
$r \in C_i$ satisfy $b(r) \equiv B \pmod{e_h}$.  Now let
$S_{i+1}$ be the collection of all such $r \in C_i$, set 
$p_{i+1} = e_h$, and let $A_{i+1} \equiv A_i \pmod{p_1\cdots p_i}$ and
$A_{i+1} \equiv B \pmod{p_{i+1}}$ by the Chinese Remainder Theorem.  
Then we will have that properties 1, 2, and 3 as above follow immediately
for this set $S_{i+1}$.\smallskip

If there exists a prime $p > L(1,x)$ which divides
more than\smallskip\noindent
 $|S_{i+1}| / \sqrt{\log{x}/\log\log{x}}$ of the elements of $S_{i+1}$, then
we set $t = i+1$ and we are finished.  If not, we continue constructing 
these sets $S_j$.  We are guaranteed to eventually hit upon such a prime
$p$ since all our $r_j$'s are divisible by at least one prime 
$p > L(1, x)$ by property B.

\enddocument